\def\simplelatex{\iffalse}
\documentclass[12pt]{article}
\usepackage[francais]{babel}
\simplelatex%
\else%
\usepackage[latin1]{inputenc}
\usepackage[T1]{fontenc}
\usepackage{ae}
\usepackage{aeguill}
\fi%
\usepackage{amsmath}
\usepackage{amsfonts}
\usepackage{amssymb}
\usepackage{amsthm}
\simplelatex%
\let\mathscr=\mathcal
\else%
\usepackage{mathrsfs}
\fi%
\simplelatex%
\let\url=\texttt
\else%
\usepackage{url}
\fi%
\newtheorem{prop}{Proposition}[section]
\newtheorem{thm}[prop]{Théorème}
\newtheorem*{maincubthm}{Théorème \ref{PrincipalCubic}}
\newtheorem*{maindp4thm}{Théorème \ref{PrincipalDP4}}
\newtheorem{lem}[prop]{Lemme}
\newtheorem{cor}[prop]{Corollaire}

\newcommand{\Spec}{\mathop\mathrm{Spec}\nolimits}
\newcommand{\Gal}{\mathop\mathrm{Gal}\nolimits}
\newcommand{\Pic}{\mathop\mathrm{Pic}\nolimits}
\newcommand{\Br}{\mathop\mathrm{Br}\nolimits}
\newcommand{\card}{\mathop\mathrm{card}\nolimits}
\newcommand{\Div}{\mathop\mathrm{Div}\nolimits}
\newcommand{\canK}{\mathrm{K}}
\simplelatex%
\let\bimu=\mu
\let\mathbi=\mathbf
\else%
\DeclareFontFamily{OML}{cmmib}{\skewchar\font127 }
\DeclareFontShape{OML}{cmmib}{m}{it}%
        {<5><6><7><8><9>gen*cmmib%
         <10><10.95>cmmib10%
         <12><14.4><17.28><20.74><24.88>cmmib12%
         }{}
\DeclareSymbolFont{biletters}{OML}{cmmib}{m}{it}
\DeclareSymbolFontAlphabet{\mathbi}{biletters}
\DeclareMathSymbol{\bimu}{\mathord}{biletters}{"16}
\fi%
\simplelatex%

\fi%
\begin{document}
\title{Surfaces de del Pezzo sans point rationnel sur un corps de
dimension cohomologique un}
\author{Jean-Louis Colliot-Thélène  et David A. Madore}
\maketitle

\begin{abstract}
Pour chaque entier $d=2,3,4$, il existe un corps $F$
de dimension cohomologique $1$ et une surface de del Pezzo de degré $d$
sur $F$ sans zéro-cycle de degré $1$, en particulier sans point
rationnel. Les démonstrations utilisent le théorème de Merkur'ev et
Suslin, le théorème de Riemann-Roch sur une surface  et la formule du
degré de Rost.
\end{abstract}

{\em Mots-clés} : points rationnels, surfaces de del Pezzo, corps 
de dimension
cohomologique $1$, théorème de Riemann-Roch, formule du degré de Rost.

Mathematical Subject Classification (MSC 2000)
Primary: 14G05, 12G10, 14C35;
Secondary:  14J26, 11E76

\section{Introduction}

Soit $F$ un corps. D'après Enriques, Manin, Iskovskikh, Mori
(voir \cite{KollarBook} III.2.1), toute
$F$-surface projective, lisse, géomé\-tri\-quement connexe qui est
(géo\-mé\-tri\-quement)
\emph{rationnelle}, i.e. devient birationnelle au plan projectif après
extension finie convenable de $F$, est, sur $F$, birationnelle
à une surface de l'un des types suivants :
\begin{itemize}
\item[(1)] surface fibrée en coniques au-dessus d'une conique ;
\item[(2)] surface de del Pezzo\footnote{Pour ce qui est
des propriétés
générales des surfaces de del Pezzo, nous renvoyons le lecteur
à \cite{KollarBook} III.3 et \cite{Manin} III et IV.}
de degré $d$, avec $1 \leq d \leq 9$.
\end{itemize}

     Une inspection de cette classification
avait ainsi permis d'observer que lorsque
le corps $F$
est un corps $C_1$ au sens de Lang (cf. \cite{SerreGC} II.§3.2)
toute telle surface possède
un point rationnel (voir par exemple \cite{KollarBook} IV.6.8).
Un corps $C_1$ est de dimension cohomologique $\leq 1$
au sens de Tate et Serre (c'est-à-dire que $\Br L=0$ pour toute
extension $L/F$ finie, cf. \cite{SerreGC} II.§3).
On peut se demander si plus généralement
toute surface rationnelle sur un corps de dimension cohomologique $1$
possède un point rationnel. C'est le cas pour les surfaces fibrées
en coniques au-dessus d'une conique. C'est le cas pour les surfaces
de del Pezzo de degré au moins égal à $5$. C'est aussi le cas
pour le degré $1$, car sur tout corps toute telle surface possède un
point rationnel évident (le point base du système anticanonique).

Le but de cet article est d'établir qu'il n'en est pas de même
en les autres degrés :

\begin{thm}
Il existe une surface cubique lisse $X$ sur $\mathbb{Q}$ et un
corps $F$ de dimension
cohomologique $1$ et de caractéristique zéro, tels que $X$ n'ait
pas de point rationnel sur $F$.

Plus précisément, $X$ ne
possède de point dans aucune extension de
degré premier à $3$
de $F$.
\label{PrincipalCubic}
\end{thm}

\begin{thm}
Il existe une intersection complète lisse $X$
de deux quadriques  dans $\mathbb{P}^4_{\mathbb{Q}}$ et un corps $F$
de dimension
cohomologique $1$ et de caractéristique zéro, tels que $X$ n'ait
pas de point rationnel sur $F$.

Plus précisément, $X$ ne
possède de point dans aucune extension de
degré impair de $F$.
\label{PrincipalDP4}
\end{thm}

Partant de ce dernier théorème, il est facile de donner
un exemple de surface de del Pezzo de degré $2$
sur $F$ sans point rationnel. Pour une autre approche,
on se reportera au  paragraphe 5 ci-dessous.

Le théorème \ref{PrincipalDP4} fournit un contre-exemple
au théorème 2 de \cite{Brumer}, dont il
était connu depuis 1981 que la démonstration est
incorrecte\footnote{L'erreur se trouve dans le dernier lemme,
p. 681 de l'article : un élément du groupe de Brauer
du corps des fractions de la droite projective qui est
non ramifié en toute place sauf peut-être en une place
de degré $2$ peut ne pas être constant.}.

Sur un corps parfait de dimension cohomologique $1$, tout
espace homogène non vide (non nécessairement principal) d'un
groupe linéaire connexe possède un point rationnel
(\cite{SerreGC}, III.2.4).
On a donc :
\begin{cor}
Pour chacune des surfaces $X$ des théorèmes
\ref{PrincipalCubic} et \ref{PrincipalDP4}, il n'existe aucune
application rationnelle définie\footnote{Pour les surfaces cubiques,
la démonstration donne encore ceci sur un corps $p$-adique.}
sur $\mathbb{Q}$ d'un espace homogène
d'un $\mathbb{Q}$-groupe linéaire connexe vers $X$.
\end{cor}
Ceci ôte l'espoir immédiat de ramener l'arithmétique de telles
surfaces à celle des groupes algébriques linéaires.

\medskip

Soit $K$ un corps algébriquement clos. Une $K$-variété
séparablement rationnellement connexe est une $K$-variété
pour laquelle il existe un $K$-mor\-phisme $f\colon\mathbb{P}^1_K \to X$
tel que l'image inverse $f^*T_X$ du fibré tangent $T_X$
soit un fibré ample sur $\mathbb{P}^1_K$.
En dimension un, on trouve la droite projective.
En dimension deux, la classe de ces variétés coïncide
avec la classe des surfaces (projectives lisses) rationnelles.
Soit $F$ un corps quelconque. Une $F$-variété projective lisse
géométriquement connexe est dite séparablement rationnellement
connexe si elle l'est après passage à une clôture algébrique de
$F$.

Les théorèmes \ref{PrincipalCubic} et \ref{PrincipalDP4} sont à
mettre en perspective avec plusieurs résultats récents.

Soit $F=k(C)$ le corps des fonctions d'une courbe projective et lisse
sur un corps $k$ algébriquement clos.  Graber, Harris et Starr
(\cite{GraberHarrisStarr}, en caractéristique zéro) et de Jong et
Starr (\cite{DeJongStarr}, en toute caractéristique) ont établi que
toute $F$-variété projective, lisse, (géométriquement)
séparablement rationnellement connexe possède un point rationnel.

Soit $F=\mathbb{F}$ un corps fini. H. Esnault a établi
dans \cite{Esnault} que toute $\mathbb{F}$-variété projective, lisse,
(géométriquement)  rationnellement connexe par chaînes (sur
un  corps assez grand)
possède un
point $\mathbb{F}$-rationnel.

C'est une question ouverte de savoir si
ces résultats sont chacun un cas particulier d'un
énoncé sur les corps $C_1$.

Rappelons que tout corps $C_1$ est de dimension cohomologique  $\leq 1$.
La réciproque n'est pas vraie, comme l'a montré
J. Ax. Dans \cite{Ax} il a donné
un exemple d'hypersurface quintique lisse
dans $\mathbb{P}^9_F$, sans point rationnel, sur un corps $F$
de caractéristique zéro de dimension cohomologique 1.
   Une telle hypersurface est
une variété de Fano,
donc rationnellement connexe (\cite{KollarBook},  V.2.13).

Un zéro-cycle  sur une $F$-variété
$X$ est une combinaison
linéaire à coefficients entiers
de points fermés de $X$. Le
degré d'un tel zéro-cycle
$\sum_in_iP_i$ est l'entier
$\sum_in_i[F(P_i):F] \in \mathbb{Z}$.

L'hypersurface quintique d'Ax possède (de façon évidente) un
zéro-cycle  de degré $1$ (il en est de même des autres
hypersurfaces construites par Ax dans le même article, et qui
établissent que pour $r>0$, un corps peut être de dimension
cohomologique 1 sans être $C_r$).

   Kato et Kuzumaki \cite{Kato}
ont posé la question de savoir si,
sur un corps $F$ de dimension cohomologique $1$,
    toute hypersurface dans $\mathbb{P}^n_F$ de degré $d\leq n$
possède un zéro-cycle de degré $1$.
     Le théorème \ref{PrincipalCubic}
montre qu'il n'en est rien\footnote{L'article \cite{Merk-uinvar}
de Merkur'ev, sur le $u$-invariant des formes quadratiques, fournit
déjà une réponse négative à une question plus générale
de \cite{Kato}.}.

\medskip

On comparera les théorèmes \ref{PrincipalCubic} et \ref{PrincipalDP4}
avec le résultat suivant (\cite{CTThmNinety}, théo\-rème A (iv)) :
pour
$X$ une surface (géométriquement) rationnelle sur un corps $k$ de
caractéristique $0$ et de dimension cohomologique $\leq 1$, tout
zéro-cycle de degré $0$ est rationnellement équivalent à $0$.

\smallskip

Les théorèmes fournissent des exemples de surfaces rationnelles $X$
sur un corps
de dimension cohomologique $1$ telles que $H^1_{\mathrm{Zar}} (X,
\mathscr{H}^3(\bimu_p^{\otimes 2})) \neq 0$, où $\bimu_p$ est le
groupe des racines $p$-ièmes de l'unité et
$\mathscr{H}^3(\bimu_p^{\otimes 2})$ est le faisceau de Zariski
associé au préfaisceau $U\mapsto H^3_{\mathrm{\acute et}} (U,
\bimu_p^{\otimes 2})$ (pour $p=3$, resp. $p=2$).

\medskip

Dans \cite{Ducros}, Ducros montre que pour toute courbe
de genre $0$ sans point rationnel sur un corps $k$ il
existe un corps $F$, contenant $k$, de dimension cohomologique $2$,
sur lequel la courbe n'a pas de point. Il utilise pour cela les
théorèmes de Merkur'ev et Suslin (\cite{MerkurjevSuslin}).
Nous utilisons aussi  les résultats de Merkur'ev
et Suslin,
mais nous avons également  recours à d'autres ingrédients.
Pour les surfaces cubiques, nous utilisons la cohomologie
galoisienne et faisons appel
à la formule du degré de Rost (cf. \cite{MerkurjevLens}).
Pour les surfaces de del Pezzo de degré $4$, nous
utilisons le théorème de Riemann-Roch pour les faisceaux inversibles
et un théorème d'Amer (\cite{Amer}) et Brumer (\cite{Brumer}) ---
théorème \ref{BrumerAmerTheorem} ci-dessous. Les paragraphes 3 et 4
de cette note (constituant respectivement la démonstration des
théorèmes
\ref{PrincipalCubic} et \ref{PrincipalDP4}) sont indépendants.
Au paragraphe 5, nous donnons une variante de la démonstration
du théorème \ref{PrincipalDP4}) qui nous a été communiquée
par J. Kollár, et s'applique à d'autres situations ;
on utilise ici le théorème de Riemann-Roch pour les fibrés
vectoriels sur une courbe.

\bigskip

Nous remercions vivement D. Coray, P. Gille,  A. Merkur'ev,
J-P. Serre, B. Totaro et tout particulièrement  J. Kollár
pour diverses remarques sur des versions antérieures de ce texte.

\section{Construction
de corps de dimension cohomologique 1}

Pour construire un corps de dimension cohomologique $1$, la technique
utilisée consiste, pour un certain nombre premier $p$ (ici, $p$ vaut
$2$ ou $3$ selon le cas), à partir d'un corps $k_0$ quelconque et à
ajouter successivement (et de façon transfinie, s'il le faut)
une extension algébrique de corps maximale de degré premier à $p$
et les corps de fonctions de
toutes les variétés de Severi-Brauer d'indice $p$.  Nous présentons
dans ce paragraphe les résultats relatifs à cette construction.

On a le résultat classique suivant (\cite{SerreGC} I.§4.1, prop. 21) :
\begin{prop}
Soient $p$ un nombre premier et $E$ un corps de caractéristique
différente de $p$.  On suppose satisfaites les deux hypothèses :
\begin{itemize}
\item[(1)] $E$ n'a pas d'extension finie de
degré premier à $p$.
\item[(2)] Le groupe de $p$-torsion de $\Br E$ est nul, ce qui
équivaut (grâce à (1))  à l'hypothèse
$H^2(E,\mathbb{Z}/p\mathbb{Z})=0$.
\end{itemize}
Alors il en est de même de celui de
toute extension finie de $E$,
et la dimension cohomologique de $E$ est au plus $1$.
\label{GaloisCohomologyFact}
\end{prop}

L'énoncé suivant est un peu moins classique :
\begin{thm} Soient $p$ un nombre premier et $E$ un corps de
caractéristique différente de $p$. Si $E$ ne satisfait pas les
hypothèses de la proposition précédente, alors soit $E$ possède
une extension finie de corps de degré premier à $p$, soit $E$
contient une racine primitive
$p$-ième $\zeta_p$ de $1$,  et il existe une algèbre cyclique
$(a,b)_{\zeta_p}$ non triviale sur $E$, ce qui définit une
variété de Severi-Brauer non triviale d'indice $p$ sur $E$.
\label{OtherGaloisCohomologyFact}
\end{thm}

Cet énoncé résulte du théorème de Merkur'ev-Suslin,
qui implique que pour $E$ contenant $\zeta_p$, le groupe
$H^2(E,\mathbb{Z}/p\mathbb{Z})$ est engendré par les cup-produits
d'éléments de $H^1(E,\mathbb{Z}/p\mathbb{Z})$.
Le fait ci-dessus,  moins difficile,
est établi à un cran intermédiaire de
la démonstration de Merkur'ev-Suslin (cf. \cite{MerkurjevSuslin},
(10.6) et (10.7)).

Le corps de dimension cohomologique $1$ sera donné par le résultat
suivant, qui découle des deux précédents :
\begin{thm}
Soient $p$ un nombre premier, $k_0$ un corps de caractéristique
différente de $p$, et $\mathbi{P}$ une propriété, stable par
isomorphisme, des extensions de
$k_0$.  On suppose que
\begin{itemize}
\item[(0)] la propriété $\mathbi{P}(k_0)$ est satisfaite,
\item[(1)] si $\mathbi{P}(k)$ est satisfaite, alors $\mathbi{P}(k')$
      est satisfaite pour toute extension algébrique $k'$ de $k$ de degré
      premier à $p$,
\item[(2)] si $\mathbi{P}(k)$ est satisfaite, alors, pour
toute variété de Severi-Brauer $Y$ associée à une algèbre
simple centrale
      d'indice $p$ sur $k$, la propriété $\mathbi{P}(k(Y))$
      est satisfaite pour  le corps des fonctions $k(Y)$ de $Y$.
\item[(3)] si $\mathbi{P}(k_\iota)$ est satisfaite pour tout membre
    $k_\iota$ d'un système inductif filtrant $(k_\iota)_{\iota\in I}$,
    alors $\mathbi{P}(\varinjlim k_\iota)$ est satisfaite.
\end{itemize}
Dans ces conditions, il existe un corps $F$, extension de $k_0$, de
dimension cohomologique $\leq 1$,
de groupe de Galois absolu un
pro-$p$-groupe,
tel que la propriété
$\mathbi{P}(F)$ soit satisfaite.
\label{InductionTheorem}
\end{thm}
Ce théorème est assez connu.  Cependant, faute de démonstration
satisfaisante dans le cadre général, nous en donnons une ci-dessous
pour la commodité du lecteur.
\begin{proof}[Démonstration du théorème \ref{InductionTheorem}]
Partant du corps $k_0$, on construit par induction transfinie sur
$\alpha$ parcourant la classe des ordinaux (en fait on va donner
ci-dessous une borne supérieure sur les $\alpha$ qui nous intéressent)
des corps $k_\alpha$, extensions de $k_0$, vérifiant
$\mathbi{P}(k_\alpha)$ pour tout $\alpha$, et avec $k_0 \subseteq
k_\beta \subseteq k_\alpha$ lorsque $\beta\leq\alpha$.

Lorsque $\delta$ est un ordinal limite, on définit $k_\delta$ comme la
limite inductive (ou, abusivement, la réunion) des $k_\alpha$ pour
$\alpha<\delta$ : comme chacun des $k_\alpha$ vérifie
$\mathbi{P}(k_\alpha)$, la propriété $\mathbi{P}(k_\delta)$ est
également vérifiée (d'après (3) ci-dessus).  Reste à expliquer
comment on construit $k_{\alpha+1}$ à partir de $k_\alpha$.

Soit $\alpha$ un ordinal pour lequel $k_\alpha$ a été construit.
Si $k_\alpha$  vérifie les hypothèses de la
proposition \ref{GaloisCohomologyFact},
c'est-à-dire
s'il n'existe pas d'extension $k'/k_\alpha$ de degré premier à $p$ et
que $H^2(k_\alpha,\mathbb{Z}/p\mathbb{Z}) = 0$, alors la dimension
cohomologique de
$k_\alpha$ est $\leq 1$ d'après la proposition en question : on pose
$k_{\alpha+1} = k_\alpha$ ;  dans ce cas $F = k_\alpha$ et on a fini.
Si ces hypothèses ne sont pas satisfaites, d'après le
théorème \ref{OtherGaloisCohomologyFact}, on peut trouver soit une
extension $k'/k_\alpha$ finie de degré premier à $p$, soit une
variété de Severi-Brauer $Y$ d'indice $p$ sur $k_\alpha$  associée
à un symbole $(a,b)_\zeta$ (avec $\zeta$ racine $p$-ième de l'unité
dans
$k_\alpha$) :
    on choisit pour $k_{\alpha+1}$ l'un
ou l'autre corps $k'$ ou $k(Y)$ (indifféremment, si les deux
conviennent), et les propriétés (1) et (2) ci-dessus assurent que
$\mathbi{P}(k_{\alpha+1})$ est encore satisfaite.

Il reste maintenant à expliquer pourquoi cette construction transfinie
se termine avant un ordinal que nous allons expliciter ($\kappa$
ci-dessous), c'est-à-dire qu'on ne construit pas une classe propre de
corps de plus en plus gros, mais qu'au contraire on a bien $k_\alpha =
k_\beta$ pour $\beta\leq\alpha$ lorsque $\beta$ est assez grand (donc
$k_\alpha$ de dimension cohomologique $\leq 1$ par construction).
Remarquons que le foncteur $H^2(\cdot,\mathbb{Z}/p\mathbb{Z})$
commute aux limites inductives
filtrantes : ainsi, lorsque $\delta$ est un ordinal limite,
$H^2(k_\delta,\mathbb{Z}/p\mathbb{Z})$
    est bien la limite des $H^2(k_\alpha,\mathbb{Z}/p\mathbb{Z})$ pour
$\alpha<\delta$.

Soit $\kappa$ le cardinal successeur de $\card k_0 + \aleph_0$ (c'est
donc un cardinal régulier indénombrable).  Effectuons une fois
pour toutes un choix d'une clôture algébrique $\bar k_\alpha$
de $k_\alpha$ pour chaque $\alpha$, de façon compatible aux inclusions
$k_\alpha\subseteq k_\beta$.
Appelons $E_\alpha$ l'ensemble produit cartésien d'ensembles
$(H^2(k_\alpha,\mathbb{Z}/p\mathbb{Z}))\times (\bar k_\alpha/k_\alpha)$
(où $\bar k_\alpha/k_\alpha$ est le $k_\alpha$-espace vectoriel quotient
de $\bar k_\alpha$ par $k_\alpha$), et
$\ast$ l'élément $(0,0)$ de ce produit.  On voit alors (par
induction sur $\alpha$) que $\card E_\alpha < \kappa$ pour tout
$\alpha<\kappa$ (puisque c'est le cas évidemment
pour $H^2( k_\alpha,\mathbb{Z}/p\mathbb{Z})$
et pour $\bar k_\alpha$).  Supposons par l'absurde que $k_\alpha$
n'est jamais de dimension cohomologique $\leq 1$ pour $\alpha<\kappa$,
c'est-à-dire que $k_{\alpha+1}$ est toujours une extension propre de
$k_\alpha$ ; alors il existe un élément différent de $\ast$ dans
$E_\alpha$ qui s'envoie sur $\ast$ dans $E_{\alpha+1}$.  Une
contradiction résulte directement du
lemme \ref{CombinatorialLemma} ci-dessous.
\end{proof}

\begin{lem}
Soit $\kappa$ un cardinal régulier indénombrable, et pour chaque
$\alpha<\kappa$ soit $E_\alpha$ un ensemble pointé par un élément
$\ast\in E_\alpha$, et $\varphi_\alpha\colon E_{\alpha}\to
E_{\alpha+1}$ tel que $\varphi_\alpha(\ast)=\ast$, et tel que pour
$\delta$ ordinal limite $E_\delta$ soit la limite du système inductif
des $E_\alpha$ pour $\alpha<\delta$ (comme ensembles pointés, avec les
$\varphi_\alpha$ pour morphismes).  On suppose de plus que (1) $\card
E_\alpha<\kappa$ pour tout $\alpha<\kappa$, et que
(2) pour tout $\alpha<\kappa$ il existe $x\in E_\alpha$
différent de $\ast$ tel que $\varphi_\alpha(x)=\ast$.
Alors il y a une contradiction.
\label{CombinatorialLemma}
\end{lem}
\begin{proof}
Appelons $N_\alpha\subseteq E_\alpha$ l'ensemble des
éléments de $E_\alpha$ qui s'envoient sur $\ast$ dans la limite
inductive $E_\kappa$ de tous les $E_\alpha$ avec $\alpha<\kappa$.
Manifestement les $N_\alpha$ vérifient les mêmes hypothèses que les
$E_\alpha$.  Soit maintenant $\gamma_0=0$, et par récurrence sur le
naturel $k$, soit $\gamma_{k+1}$ le plus petit ordinal $\alpha<\kappa$
tel que tous les éléments de $N_{\gamma_k}$ aient pour image $\ast$
dans $N_\alpha$ (un tel $\alpha$ existe puisque chaque élément donné
de $N_{\gamma_k}$ s'envoie sur $\ast$ dans $N_\kappa$, donc aussi dans
un $N_\beta$ avec $\beta<\kappa$, et en prenant la borne supérieure
de ces $\beta$-là, qui porte sur un ensemble de cardinal $<\kappa$, on
obtient le $\gamma_{k+1}<\kappa$ recherché, l'inégalité stricte
provenant de
ce que $\kappa$ est régulier).  Soit maintenant $\gamma$ la borne
supérieure des $\gamma_k$ pour $k\in\mathbb{N}$ : on a $\gamma<\kappa$
car $\kappa$ est indénombrable régulier.  Or si $x$ est un élément de
$N_\gamma$, il appartient à un $N_{\gamma_k}$ pour $k\in\mathbb{N}$,
donc s'envoie sur $\ast$ dans $N_{\gamma_{k+1}}$, donc dans
$N_\gamma$, ce qui prouve que $N_\gamma=\{\ast\}$.  Mais ceci
constitue une contradiction, car il existe par hypothèse un
$x\neq\ast$ dans $E_\gamma$ qui s'envoie sur $\ast$ dans
$E_{\gamma+1}$, donc $x\in N_\gamma$, contredisant $N_\gamma=\{\ast\}$.
\end{proof}

\textbf{Remarques :} Pour résumer cet épisode combinatoire, disons que
si on part du corps $k=k_0$ et qu'on itère de façon transfinie
l'opération de passer à une extension finie de degré premier à
$p$ ou celle de passer au corps des fonctions d'une variété de
Severi-Brauer, tant que l'une d'elles au moins est possible, n'importe
quel choix de ces extensions finira par donner, en moins de $\kappa$
étapes, un corps $F$ comme souhaité.  Notons que si on se contente de
prendre des extensions de degré premier à $p$, on obtient
le corps fixe d'un pro-$p$-Sylow dans le groupe de Galois : on peut
donc présenter différemment la construction transfinie (c'est sans
doute plus classique et évite d'utiliser la notion d'ordinaux, mais
l'argument semble plus \textit{ad hoc} donc moins agréable) en effectuant
alternativement l'opération de prendre le corps fixe d'un
pro-$p$-Sylow et celle de passer au corps des fonctions de toutes les
variétés de Severi-Brauer (à isomorphisme près)
sur le corps précédent (comparer avec la
construction utilisée dans \cite{MerkurjevSuslin} (11.4)
et avec celle utilisée dans \cite{Ducros}).  Notons que le corps
$F$ finalement obtenu a un cardinal inférieur à $\kappa$,
c'est-à-dire, au
plus égal au cardinal de $k$ si $k$ est infini, au plus dénombrable
sinon ; avec l'utilisation que nous ferons ci-dessous de ce théorème,
$F$ est de cardinal $\leq 2^{\aleph_0}$,
mais un argument simple  permet de voir qu'un sous-corps dénombrable suffit
déjà.

\bigskip

Nous terminons ce paragraphe avec un résultat qui, s'il n'est pas
directement utile dans ce qui suit, sert au moins à motiver certaines
conditions introduites :
\begin{prop}
Soient $k$ un corps de caractéristique $0$ et $X$ une variété sur $k$
projective, lisse, géométriquement connexe, telle que $H^1(X,
\mathcal{O}_X) = 0$ (hypothèse satisfaite, par exemple, si $X$ est
$\bar k$-rationnelle).  Alors pour toute extension $K$ de $k$ dans
laquelle $k$ est algébriquement fermé, la flèche naturelle $(\Pic
X_{\bar k})^{\Gal(\bar k/k)} \to (\Pic X_{\bar K})^{\Gal(\bar K/K)}$
est un isomorphisme.
\label{PicardProperty}
\end{prop}
\begin{proof}
Pour alléger les notations, soit $G_k = \Gal(\bar k/k)$ le groupe de
Galois absolu de $k$, et $G_K = \Gal(\bar K/K)$ celui de $K$, et soit
de plus $H = \Gal(\bar K/\bar kK)$, de sorte que $G_k = \Gal(\bar kK /
K) = G_K/H$.  Sous l'hypothèse $H^1(X, \mathcal{O}_X) = 0$, on sait
(cf. \cite{Hartshorne}, chapitre III exercice 12.6 p. 192) que les flèches
naturelles $\Pic X_{\bar k} \to \Pic X_{\bar k K} \to \Pic X_{\bar K}$
sont des isomorphismes. Ainsi
$H$ agit trivialement sur $\Pic X_{\bar K}$.
Ainsi, $(\Pic X_{\bar k})^{G_k}\to (\Pic X_{\bar K})^{G_K}$ est un
isomorphisme.
\end{proof}

\section{Surfaces de del Pezzo de degré $3$}

Pour démontrer le théorème \ref{PrincipalCubic}, nous aurons besoin
de certains faits connus sur les surfaces cubiques, qui sont résumés
dans le théorème suivant, dû pour l'essentiel à B. Segre :
\begin{thm}
Soit $X$ une surface cubique lisse sur un corps $k$ parfait.  Soit
$\bar k$ une clôture algébrique de $k$, soit $\bar X = X\times_{\Spec
k} \Spec\bar k$ et soit $G = \Gal(\bar k/k)$ le groupe de Galois de $\bar
k$ sur $k$.  Alors les deux affirmations suivantes sont équivalentes :
\begin{itemize}
\item[(a)] La surface $X$ est $k$-minimale.
\item[(b)] Le groupe $(\Pic\bar X)^G$ des invariants par $G$ du groupe
      de Picard $\Pic\bar X$ de $\bar X$ est (libre) de rang 1.

\end{itemize}
Lorsqu'elles sont satisfaites, l'injection naturelle $\Pic X
\hookrightarrow (\Pic\bar X)^G$ est une bijection,
ces deux groupes sont engendrés par la classe du fibré
canonique $\canK_X$ (classe dont l'opposée est induite par le
plongement naturel de
la surface cubique $X $ dans $ \mathbb{P}^3_k$),  et
l'application naturelle $\Br k\to \Br
k(X)$ est injective.  \par De plus, si $X$ est diagonale, c'est-à-dire
donnée par une équation de la forme $a_0 T_0^3 + a_1 T_1^3 + a_2 T_2^3
+ a_3 T_3^3 = 0$, alors ces conditions sont satisfaites dès que la
suivante l'est :
\begin{itemize}
\item[(ii)] aucun des  rapports $a_m
a_n/a_p a_q$, avec $(m,n,p,q)$ permutation de
$(0,1,\penalty0 2,3)$, n'est dans $k^{\times 3}$.
\end{itemize}
\label{SegreWonderfulTheorem}
\end{thm}
\begin{proof}
L'équivalence de (a) et (b) est un théorème de B. Segre
(voir \cite{Manin}, théorème IV.28.1, page
151). On a les inclusions naturelles de groupes abéliens libres
de rang un $\mathbb{Z} \canK_X
\hookrightarrow
\Pic(X) \hookrightarrow (\Pic\bar X)^G$, et la classe
de $\canK_X$ n'est pas divisible dans $\Pic\bar X$
(son nombre d'intersection avec une des droites tracées sur
$\bar X$ vaut $-1$) : ces inclusions sont donc des égalités.

Pour toute $k$-variété projective, lisse, géométriquement intègre
$X$, de corps des fractions $k(X)$, la cohomologie galoisienne
de la suite exacte naturelle
$1 \to \bar k(X)^\times/\bar k^\times \to \Div \bar X \to \Pic \bar X \to
0$ donne la  suite exacte classique
\[
0\to \Pic X \to \Pic\bar X^G \to \Br k \to \Br k(X)
\]
Ceci établit l'injectivité de la flèche $\Br k\to \Br k(X)$ lorsque
les conditions (a) et (b) sont satisfaites.

Enfin, le fait que la condition (ii) garantisse
la $k$-minimalité de $X$
est de nouveau un résultat de B. Segre (voir \cite{Manin} 
exercice III.21.10).
\end{proof}

Le fait suivant, qui énonce un résultat tout à fait parallèle au
théorème \ref{Cub2Iterating} plus bas, pour une autre sorte
d'extension (cf. (1) du théorème \ref{InductionTheorem}),
est évident :
\begin{lem}
Soient $k$ un corps de caractéristique $0$ et $X$ la surface cubique
lisse dans $\mathbb{P}^3_k$ définie par l'équation $a_0 T_0^3 + a_1
T_1^3 + a_2 T_2^3 + a_3 T_3^3 = 0$, où $a_0,a_1,a_2,a_3 \in k^\times$.
On suppose que
\begin{itemize}
\item[(i)] $X$ n'admet de point dans aucune extension de corps
$k'/k$ de degré premier à $3$,
\item[(ii)] aucun des  rapports $a_m
a_n/a_p a_q$, avec $(m,n,p,q)$ permutation de
$(0,1,\penalty0 2,3)$, n'est dans $k^{\times 3}$.
\end{itemize}
Alors, si $k'$ est
une extension finie de $k$ de degré non
multiple de $3$, la surface $X$ étendue au corps $k'$ vérifie toujours
les propriétés (i) et (ii).
\label{Cub1Iterating}
\end{lem}

En combinant ces divers résultats et la formule du degré de Rost,
nous établissons maintenant le résultat technique principal (qui vise
le (2) du théorème \ref{InductionTheorem}).

\begin{thm}
Soient, comme dans le lemme \ref{Cub1Iterating}, $k$ un
corps de carac\-téris\-tique $0$ et $X$ la surface cubique lisse dans
$\mathbb{P}^3_k$ définie par l'équation $a_0 T_0^3 + a_1 T_1^3 + a_2
T_2^3 + a_3 T_3^3 = 0$, où $a_0,a_1,a_2,a_3 \in k^\times$.  Supposons
les deux conditions suivantes satisfaites
\begin{itemize}
\item[(i)] $X$ n'admet de point dans aucune extension de corps $k'/k$ de
degré premier à $3$,
\item[(ii)] aucun des  rapports $a_m
a_n/a_p a_q$, avec $(m,n,p,q)$ permutation de
$(0,1,\penalty0 2,3)$, n'est dans $k^{\times 3}$.
\end{itemize}

     Soient $D$ une algèbre simple centrale sur $k$
d'indice $3$, $Y$ la surface de Severi-Brauer
associée à $D$ et $k(Y)$ son corps des fonctions.

Alors
la $k(Y)$-surface $X\times_{\Spec k}\Spec k(Y)$
vérifie encore les propriétés (i) et (ii).

(En particulier, elle n'a pas de point sur $k(Y)$, c'est-à-dire qu'il
n'existe pas d'application rationnelle $Y\dasharrow X$ définie sur $k$.)
\label{Cub2Iterating}
\end{thm}

\begin{proof}
Nous utiliserons les notations et les résultats du rapport de
Merkur'ev \cite{MerkurjevLens} sur la formule du degré de Rost, que
nous rappelons brièvement ici ; voir aussi
\cite{MerkurjevSteenrod}, \cite{RostDegreeOther} et \cite{RostDegree}.

A toute variété $Z$  sur un corps $k$ on  attache
un entier naturel $n_Z$, appelé l'indice de $Z$ sur $k$ : c'est par
définition le plus grand commun diviseur des degrés, sur
$k$,  des corps résiduels
$k(P)$ des points fermés $P$ de $Z$ ; c'est aussi le p.g.c.d. des
degrés des extensions finies $K/k$ telles que $X(K)\neq \varnothing$.
De plus, pour $p$ un nombre
premier et
$Z$ projective, on définit un invariant $\eta_p(Z)\in
\mathbb{Z}/n_Z\mathbb{Z}$ (tué par $p$).  Cet invariant vérifie la
formule du degré (\cite{MerkurjevLens}, théorème 4.1) : si
$Y$ et $X$ sont des $k$-variétés projectives intègres de même
dimension, et si  $f\colon Y\to X$ est un $k$-morphisme et $\deg f$
désigne le degré de $f$ (qui vaut $0$ si $f$ n'est pas dominant, et
$[k(Y):k(X)]$ si $f$ l'est), alors
$n_X$ divise $n_Y$ et $\eta_p(Y) = \deg f\cdot \eta_p(X) \in
\mathbb{Z}/n_X\mathbb{Z}$.  Par ailleurs, on sait calculer $n_Y = p^n$
et $\eta_p(Y)=p^{n-1}\in \mathbb{Z}/p^n\mathbb{Z}$ lorsque $Y$ est la
variété de Severi-Brauer d'une algèbre simple centrale d'indice
$p^n$ : cf. \cite{MerkurjevSteenrod} §7.2 et \cite{MerkurjevLens}
remarque 7.5.

Dans la situation qui nous intéresse, l'hypothèse (i) faite sur $X$ se
traduit par $n_X=3$ ; par ailleurs on a $n_Y=3$ puisque l'algèbre $D$
est d'indice $3$, et de plus $\eta_3(Y)=1 \in\mathbb{Z}/3\mathbb{Z}$
(ce qui nous intéresse est que cet invariant ne soit pas nul).

Soit $K/k(Y)$ une extension finie de corps de degré premier à 3.
Soit $Y'$ la normalisation (projective, intègre)
de $Y$ dans $K$. On dispose donc
du $k$-morphisme fini de $k$-variétés intègres $p\colon Y' \to Y$.
Si  $X\times_k{k(Y)}$ possède un point dans $K$,
alors il existe une application rationnelle $h$, définie sur $k$
de $Y'$ vers $X$. Soit $Y''\subseteq Y' \times_kX$ l'adhérence du graphe de
l'application rationnelle $h$. C'est une $k$-variété projective, intègre.
On
dispose alors du $k$-morphisme composé $q\colon Y'' \to Y$,
de degré premier à 3, et du $k$-morphisme $r\colon Y''\to X$ donné
par la seconde projection.
En appliquant la formule du degré régulier
(\cite{MerkurjevLens}, théorème 4.1) au morphisme $q$,  on obtient
$\eta_p(Y'') = \deg q\cdot \eta_p(Y) \in
\mathbb{Z}/n_Y\mathbb{Z}=\mathbb{Z}/3\mathbb{Z}$.
En appliquant la même formule au morphisme $Y''\to X$,
on obtient $\eta_p(Y'') = \deg r\cdot \eta_p(X) \in
\mathbb{Z}/n_X\mathbb{Z}=\mathbb{Z}/3\mathbb{Z}$. De la première formule
on déduit  $\eta_p(Y'')\neq 0 \in \mathbb{Z}/3\mathbb{Z}$, de la seconde
on déduit alors que le degré de $r$ est premier à 3.
Ainsi l'application naturelle $r^*\colon \mathop{_3\Br}k(X)
\to \mathop{_3\Br}k(Y'')$ induite sur la $3$-torsion
des groupes de Brauer est injective.
Mais par ailleurs la condition (ii)  implique, d'après le
théorème \ref{SegreWonderfulTheorem}, que la flèche $\Br k\to \Br
k(X)$ est injective, donc la flèche naturelle $\mathop{_3\Br}k
\to \mathop{_3\Br}k(Y'')$ l'est aussi.  Mais ceci contredit le fait
que la classe de $D$ dans $\mathop{_3\Br}k$, qui n'est pas nulle,
s'envoie sur la classe nulle dans $\Br k(Y)$, et donc
aussi sur la classe nulle dans $\Br k(Y'')$.
Ceci prouve qu'il n'existe pas d'application $k$-rationnelle $f\colon
Y'\dasharrow X$ : la $k(Y)$-surface $X\times_kk(Y)$ satisfait
l'analogue de la condition (i).

Pour ce qui est de la condition (ii) pour la $k(Y)$-surface
$X\times_kk(Y)$, elle est évidente, car
$k$ est algébriquement fermé dans $k(Y)$.  On peut également appliquer
la proposition \ref{PicardProperty} prouvée plus haut de façon plus
générale : pour toute surface cubique  lisse sur $k$,
l'hypothèse (b) du  théorème \ref{SegreWonderfulTheorem} est donc
invariante par passage de $k$ à une extension $K$ dans laquelle $k$
est algébriquement fermé.
\end{proof}

Pour appliquer le théorème \ref{InductionTheorem},
nous aurons encore besoin du fait suivant (qui vise le (0) de ce
théorème) :
\begin{lem}
Il existe un corps $k$ de caractéristique $0$ et quatre
éléments  $a_0,a_1,\penalty0 a_2,a_3
\in {\mathbb{Q}}^\times$ tels que la surface cubique lisse $X$ dans
$\mathbb{P}^3_k$
définie par l'équation $a_0 T_0^3 + a_1 T_1^3 + a_2 T_2^3 + a_3 T_3^3
= 0$ vérifie les conditions (i) et (ii) du
lemme \ref{Cub1Iterating}.
\label{CubStarting}
\end{lem}

\begin{proof}
Prenons $p$ un nombre premier tel que $p\equiv 1\pmod{3}$, et soit
$\tilde a \in\mathbb{F}_p^\times$ qui ne soit pas un cube et $a
\in\mathbb{Z}_p^\times$ qui relève $\tilde a$ modulo $p$.  Posons
$k=\mathbb{Q}$ ou $k=\mathbb{Q}_p$ et considérons la surface
cubique $X$ d'équation
$T_0^3 + p T_1^3 + p^2 T_2^3 - a T_3^3 = 0$.  Le fait (ii) est
évident.  Pour montrer que $X$ n'a pas de point sur $\mathbb{Q}_p$, il
suffit de considérer les valuations des trois termes $T_0^3$, $p
T_1^3$ et $p^2 T_2^3$, qui sont distinctes (modulo $3$) : on doit
avoir $\tilde T_0^3 = \tilde a \tilde T_3^3$ après réduction,
mais comme $\tilde a$ n'est pas un cube c'est impossible.
Plus généralement, si $k'$ est une extension de $\mathbb{Q}_p$ de
degré premier à $3$, alors pour la même raison $X$ n'a pas non plus de
point sur $k'$ car l'indice de ramification est premier à $3$ (ce qui
permet de reproduire l'argument sur les valuations) et $\tilde a$ n'est
toujours pas un cube dans le corps résiduel de $k'$.
\end{proof}

(Un théorème de Coray \cite{CorayCub} affirme que, sur un corps
$p$-adique, une surface cubique $X$ ne possédant pas de point
rationnel vérifie toujours $n_X=3$, c'est-à-dire n'a de point dans
aucune extension de degré premier à $3$ --- notre hypothèse (i).
C'est une question ouverte de savoir si ce résultat vaut sur tout
corps.  Comparer avec le théorème \ref{BrumerAmerTheorem} plus bas.)

\medskip

Ces différents ingrédients nous permettent maintenant
d'établir le théorème sur les surfaces cubiques :
\begin{maincubthm}
Il existe une surface
cubique $X$ lisse dans $\mathbb{P}^3_{\mathbb{Q}}$ et un corps $F$
de caractéristique $0$ et de dimension
cohomologique $1$ tels que $X(F) =
\varnothing$.

Plus précisément, $X$ ne possède de point dans aucune
extension  de degré premier à $3$ de $F$.
\end{maincubthm}

\begin{proof}
On part d'un corps $k_0 = k$ et d'une surface cubique $X$ lisse
diagonale dans
$\mathbb{P}^3_{\mathbb{Q}}$ qui vérifie sur $ {\mathbb{Q}}$ et sur
$k$ les hypothèses (i) et (ii) du
lemme \ref{Cub1Iterating}, comme le lemme \ref{CubStarting} en
assure l'existence.

Lorsque $K$ est une extension de $k_0$, on appellera $\mathbi{P}(K)$
la propriété suivante : l'extension $X_K = X\times_{\Spec k_0} \Spec
K$ de la surface $X$ à $K$ vérifie les hypothèses (i) et (ii) du
lemme \ref{Cub1Iterating}.

Le lemme \ref{Cub1Iterating} garantit la condition (1) du
théorème \ref{InductionTheorem}.  Le théorème \ref{Cub2Iterating}
garantit la condition (2).  Le lemme \ref{CubStarting}, pris comme
point de départ, garantit la condition (0).  Enfin, la condition (3)
est tout à fait claire.  Le résultat voulu découle donc du
théorème \ref{InductionTheorem}.
\end{proof}

\section{Surfaces de del Pezzo de degré $4$}

Pour prouver le théorème \ref{PrincipalDP4},
nous utiliserons un  résultat élégant dû à Amer et Brumer
(\cite{Amer} ; \cite{Brumer}, théorème 1) :
\begin{thm}
Soit $k$ un corps (de caractéristique $>2$),
et soient $Q_1$ et $Q_2$ deux formes quadratiques
en $r$ variables sur $k$.  S'il existe une extension algébrique $k'/k$
de degré impair telle que $Q_1$ et $Q_2$ aient un vecteur isotrope
commun dans $k^{\prime r}$ alors elles en ont un déjà dans $k^r$.
\label{BrumerAmerTheorem}
\end{thm}

(Dans le cas que nous utiliserons, le résultat avait été
obtenu par Coray \cite{CorayDP4}.)

En d'autres termes,  dès qu'une intersection de
deux quadriques possède un zéro-cycle de degré impair, elle
possède un point rationnel.

Établissons maintenant que cette propriété de ne pas avoir de
zéro-cycle de degré impair est préservée par le passage au corps
des fonctions d'une conique. On comparera la démonstration du
théorème suivant avec celle du théorème \ref{Cub2Iterating},
laquelle reposait sur la formule du degré de
Rost.
\begin{thm}
Soient $k$ un corps de caractéristique $0$ et $X$ une intersection
complète lisse de deux quadriques dans $\mathbb{P}^4_k$.  On suppose
que $X(k)=\varnothing$ et que le groupe de Picard $\Pic X$ de $X$ est
réduit aux multiples de la classe $\canK_X$ du fibré canonique
(classe dont l'opposée est induite par le  plongement naturel de
la surface  $X $ dans
$\mathbb{P}^4_k$). Soit $C$
une conique lisse sur $k$.  Alors on a $X(k(C))=\varnothing$.
\label{DP4Iterating}
\end{thm}
\begin{proof}
Supposons par l'absurde qu'il existe une application ration\-nelle
$h\colon C\to X$, qui est alors un morphisme.
L'image de $h$ ne peut pas être de dimension $0$ car elle définirait
un point rationnel de $X$.  C'est donc une courbe
géométriquement intègre ${\mit\Gamma}$
définie sur $k$.  La classe $[{\mit\Gamma}]$ de ${\mit\Gamma}$ dans
$\Pic X$ s'écrit $n \canK_X$ pour un certain entier $n$ (noter que $n<0$).
D'après la formule d'adjonction (cf. \cite{Hartshorne}, chapitre V,
proposition 1.5 et \cite{SerreGACC}, chapitre IV, n°8,
proposition 5), le genre arithmétique $p_{\textrm{a}}({\mit\Gamma})$ de
${\mit\Gamma}$ est $2n(n+1)+1$ puisque l'auto-intersection de $\canK_X$ est
$4$.

La flèche $C\to{\mit\Gamma}$ se factorise par la normalisation
$\varpi\colon\tilde{\mit\Gamma}\to{\mit\Gamma}$ de
${\mit\Gamma}$.
La $k$-courbe $\tilde{\mit\Gamma}$ est géométriquement intègre.
Le genre géométrique $p_{\textrm{g}}({\mit\Gamma})$
de ${\mit\Gamma}$ est le genre de $\tilde{\mit\Gamma}$ qui est
$0$ car le genre de $C$ est $0$.  Ainsi, $p_{\textrm{g}}({\mit\Gamma})
= 0$ (qui est pair).  Le
genre arithmétique $p_{\textrm{a}}({\mit\Gamma})$ de ${\mit\Gamma}$,
comme expliqué ci-dessus, est $2n(n+1)+1$ (et il est impair).  On a la
suite exacte courte de faisceaux cohérents sur ${\mit\Gamma}$
suivante : $0\to \mathcal{O}_{\mit\Gamma} \to
\varpi_*\mathcal{O}_{\tilde{\mit\Gamma}}
\to \mathcal{F}\to 0$ (où le faisceau gratte-ciel
$\mathcal{F}$ est, par définition, le
conoyau ainsi indiqué).  De plus, d'après \cite{SerreGACC},
chapitre IV, n°7, proposition 3, le genre arithmétique de
${\mit\Gamma}$ s'écrit $p_{\textrm{a}}({\mit\Gamma}) =
p_{\textrm{g}}({\mit\Gamma}) + \delta$ où $\delta =
\sum_{Q\in{\mit\Gamma}} \delta_Q$ est la longueur
sur $k$ du faisceau
$\mathcal{F}$ : ainsi, $\delta$ est impair.  Le cycle
$\sum_{Q\in{\mit\Gamma}} \delta_Q [Q]$ (où $\delta_Q$ désigne la
longueur de la fibre en $Q$ de $\mathcal{F}$), sur ${\mit\Gamma}$,
donc sur $X$, est alors de degré impair, donc, d'après le
théorème \ref{BrumerAmerTheorem}, il y a sur $X$ un point rationnel,
une contradiction.
\end{proof}

La propriété $\Pic X = \mathbb{Z}\cdot \canK_X$ est elle-même
préservée
par passage au corps des fonctions d'une conique : ceci résulte de la
proposition \ref{PicardProperty} plus haut.  Comme pour le cas des
surfaces cubiques, on a $\Pic(X_{\bar k})^{\Gal(\bar k/k)} =
\mathbb{Z} \cdot \canK_X$ si et seulement si $\Pic X = \mathbb{Z}\cdot
\canK_X$ : on a les inclusions naturelles de groupes abéliens libres de
rang un $\mathbb{Z} \canK_X \hookrightarrow \Pic(X) \hookrightarrow
(\Pic\bar X)^{\Gal(\bar k/k)}$, et la classe de $\canK_X$
n'est pas divisible dans $\Pic\bar X$ (son nombre
d'intersection avec une  droite tracée
sur $\bar X$ vaut $-1$),
ces inclusions sont donc des égalités.

\medskip

On peut maintenant expliciter la surface considérée :
\begin{prop}
Soit $k$ un corps de caractéristique $0$.  On suppose que le degré de
$k(\sqrt{-1}, \sqrt{2}, \sqrt{3}, \sqrt{5})$ sur $k$ est $16$.
L'intersection $X\subset \mathbb{P}^4_k$ des deux quadriques
définies par les équations
\[
\begin{array}{c}
T_0^2 + T_1^2 + T_2^2 + T_3^2 + T_4^2 = 0\\
2 T_0^2 + 3 T_1^2 + 5 T_2^2 + 7 T_3^2 + 11 T_4^2 = 0
\end{array}
\]
est lisse.  Le groupe
$(\Pic X_{\bar k})^{\Gal(\bar k/k)}$ est réduit aux multiples de
la classe $\canK_X$ du fibré canonique.
\label{DP4Starting}
\end{prop}
\begin{proof}
La vérification de la lissité de $X$ est immédiate.

Dans $\bar k^5$ de coordonnées $(T_0,T_1,T_2,T_3,T_4)$, les deux
vecteurs suivants :
\[
\begin{array}{c}
(\sqrt{-6},\sqrt{10},\sqrt{-5},1,0)\\
(2\sqrt{5},-3\sqrt{-3},\sqrt{6},0,1)
\end{array}
\]
sont simultanément isotropes et orthogonaux pour les deux formes
quadratiques $T_0^2 + T_1^2 + T_2^2 + T_3^2 + T_4^2$ et $2 T_0^2 + 3
T_1^2 + 5 T_2^2 + 7 T_3^2 + 11 T_4^2$.  La droite qui relie les deux
points ainsi définis dans $\mathbb{P}^4_{\bar k}$ est donc entièrement
contenue dans $X_{\bar k}$.  Une vérification aisée permet de voir que
l'image de cette droite par les seize éléments de
$\Gal(k(\sqrt{-1}, \sqrt{2}, \sqrt{3}, \sqrt{5}) / k) \cong
(\mathbb{Z}/2\mathbb{Z})^4$ donne bien seize droites distinctes
sur $X_{\bar k}$.

On sait par ailleurs (cf. \cite{Manin}) que les classes des seize
droites tracées sur $X_{\bar k}$ engendrent le groupe de Picard.
Comme l'action de $G=\Gal(\bar k/k)$ est transitive sur celles-ci,
$(\Pic X_{\bar k})^G$ est de rang $1$.
Comme $\canK_X$ n'est pas divisible dans $\Pic
X_{\bar k}$, on a plus précisément $(\Pic X_{\bar k})^G =
\mathbb{Z}\cdot \canK_X$.
\end{proof}

\textbf{Remarque :} Nous avons donné ici un exemple numérique pour ne
pas nous embarrasser de calculs  fastidieux.  Néanmoins, si
les deux formes quadratiques considérées sont $a_0 T_0^2 + a_1 T_1^2 +
a_2 T_2^2 + a_3 T_3^2 + a_4 T_4^2$ et $b_0 T_0^2 + b_1 T_1^2 + b_2
T_2^2 + b_3 T_3^2 + b_4 T_4^2$ avec les $a_i$ et les $b_i$ tels
qu'aucun des $d_{ij} = a_i b_j - a_j b_i$ ne soit nul, alors les seize
droites sont données par
\[
\begin{array}{c}
(\varepsilon_0 \sqrt{\frac{d_{13}d_{23}d_{04}}{d_{01}d_{20}d_{34}}},
\varepsilon_1 \sqrt{\frac{d_{23}d_{03}d_{14}}{d_{12}d_{01}d_{34}}},
\varepsilon_2 \sqrt{\frac{d_{03}d_{13}d_{24}}{d_{20}d_{12}d_{34}}},
1,0)\\
(\delta \varepsilon_0 \sqrt{\frac{d_{14}d_{24}d_{03}}{d_{01}d_{20}d_{43}}},
\delta \varepsilon_1 \sqrt{\frac{d_{24}d_{04}d_{13}}{d_{12}d_{01}d_{43}}},
\delta \varepsilon_2 \sqrt{\frac{d_{04}d_{14}d_{23}}{d_{20}d_{12}d_{43}}},
0,1)\\
\end{array}
\]
(c'est-à-dire que les deux vecteurs en question dans $\bar k^5$ sont
simultanément isotropes et orthogonaux pour les deux formes
quadratiques)
en parcourant les choix de signes $(\varepsilon_0, \varepsilon_1,
\varepsilon_2, \delta) \in \{\pm 1\}^4$.  Lorsque le groupe de Galois
définit une action transitive sur ces choix de signes, la même
démonstration que ci-dessus s'applique.

\medskip

Ces différents ingrédients nous permettent maintenant
d'établir le théo\-rème principal sur les surfaces de del Pezzo
de degré $4$ :
\begin{maindp4thm}
Il existe une surface $X$ intersection
complète lisse de deux quadriques dans $\mathbb{P}^4_\mathbb{Q}$ et
un corps $F$ de caractéristique $0$ et de dimension
cohomologique $1$ tels que
$X(F) = \varnothing$.

Plus précisément, $X$ ne possède de
point dans aucune extension
de degré impair de $F$.
\end{maindp4thm}

\begin{proof}
On part de la surface $X$ donnée par les équations de la
proposition \ref{DP4Starting} sur le corps $k_0 = \mathbb{Q}$.
Alors  $[k_0(\sqrt{-1}, \sqrt{2}, \sqrt{3},
\sqrt{5}) : k_0] = 16$ et $X(k_0) = \varnothing$.

Lorsque $K$ est une extension de $k_0$, on appellera $\mathbi{P}(K)$
la conjonction des deux propriétés suivantes :
\begin{itemize}
\item[(i)] $X_{K}(K) = \varnothing$,
\item[(ii)] $[K(\sqrt{-1}, \sqrt{2}, \sqrt{3},
\sqrt{5}) : K] = 16$.
\end{itemize}
Notons d'ores et déjà qu'en vertu de la
proposition \ref{DP4Starting}, la propriété (ii) permet de conclure
que $\Pic X_{K} = \mathbb{Z}\cdot \canK_X$.

On veut maintenant appliquer le théorème \ref{InductionTheorem}.  La
propriété (i) est préservée par extension algébrique de
degré impair
en vertu du théorème \ref{BrumerAmerTheorem} ; et la propriété (ii)
l'est aussi, de façon évidente : ceci garantit la condition (1) du
théorème \ref{InductionTheorem}.  La propriété (i) est
préservée, en
présence de la propriété (ii), par passage au corps des fonctions
d'une conique, d'après le théorème \ref{DP4Iterating} ; et la
propriété (ii) l'est aussi de façon évidente (comparer avec la
proposition \ref{PicardProperty}) : ceci garantit la condition (2).
La condition (0) est assurée dans notre choix initial.  Enfin, la
condition (3) est claire.  Le théorème \ref{InductionTheorem} permet
alors de conclure.
\end{proof}

\section{Un raffinement, et quelques applications}

J. Kollár nous a indiqué un raffinement de la démonstration
du théorème \ref{DP4Iterating}. Avec son accord, nous
décrivons son argument. L'idée supplémentaire est d'utiliser le
théorème de Riemann-Roch pour les fibrés vectoriels sur une courbe.
Cela permet de faire l'économie du  théorème \ref{BrumerAmerTheorem}.

Rappelons qu'on a défini au début de la démonstration du
théorème \ref{Cub2Iterating} l'indice $n_Z$ d'une variété $Z$ sur un
corps $k$ comme le p.g.c.d. des degrés des extensions finies $K/k$
telles que $X(K)\neq \varnothing$.

\begin{thm}
Soient $k$ un corps de caractérisitique $0$ et $X$
une $k$-surface projective, lisse, géométriquement intègre.
Soit $\canK_X \in \Pic X$ la classe canonique.
Soit $\ell$ un nombre premier qui divise l'indice $n_X$ de $X$.
Soit $C/k$ une $k$-courbe projective, lisse,
intègre, dont la caractéristique d'Euler-Poincaré
$\chi_C(\mathcal{O}_C)$ (relative au corps de base $k$) est première
à $\ell$. Soient
$D$ une
$k$-courbe lisse intègre,
$f\colon D \to C$ un $k$-morphisme fini de degré premier à $\ell$,
et $g\colon D \to X$ un $k$-morphisme. Alors l'image de $D$
est une $k$-courbe intègre ${\mit\Gamma}\subset X$ telle que $\ell$ ne
divise pas l'entier
$({\mit\Gamma}\cdot({\mit\Gamma}+\canK_X))/2$.
\label{RiemannRoch}
\end{thm}

\begin{proof}
Si $C$ possèdait un point fermé de degré (sur $k$)
non divisible par $\ell$, comme le degré de $f$ n'est pas
divisible par $\ell$, il en serait de même de $D$,
puis de $X$, ce qui est exclu. Ainsi $\ell$ divise $n_C$.

   Pour les caractéristiques
d'Euler-Poincaré cohérentes (où toutes les dimensions sont prises
sur le corps de base $k$), on a, puisque
$f$ est fini, les égalités
$\chi_D(\mathcal{O}_D)=\chi_C(f_*\mathcal{O}_D)$. Comme le morphisme $f$
est fini et plat, le faisceau cohérent $f_*\mathcal{O}_D$ est un fibré
vectoriel sur $C$ de rang $r$, où $r$ est le degré de $f$. Le
théorème de Riemann-Roch pour les fibrés vectoriels sur la
$k$-courbe
$C$ donne l'égalité
$\chi_C(f_*\mathcal{O}_D)=\deg(\bigwedge^{\max}(f_*\mathcal{O}_D))+
r\chi_C(\mathcal{O}_C)$. Comme le premier $\ell$ divise $n_C$, il divise
$\deg(\bigwedge^{\max}(f_*\mathcal{O}_D))$. On conclut
$\chi_D(\mathcal{O}_D) \equiv r\chi_C(\mathcal{O}_C) \pmod{\ell}$,
ce qui d'après les hypothèses implique
que $\chi_D(\mathcal{O}_D)$ est premier à $\ell$.

   Soit $k_D$ le corps extension finie de $k$
qui est la clôture intégrale de $k$ dans le corps des
fonctions de $D$.  Le degré de $k_D$ sur $k$ est premier à $\ell$,
puisqu'il en est ainsi de celui de $f$. L'image du morphisme
propre $g \colon D \to X$ ne saurait être de dimension zéro :
ce serait alors un point fermé de $X$ dont le degré sur $k$
diviserait celui de $k_D$ sur $k$, ce qui contredirait l'hypothèse
que $\ell$ divise $n_X$. Soit donc ${\mit\Gamma} \subset X$ la courbe
image de $g$. Soit $\varpi \colon \tilde{\mit\Gamma}\to {\mit\Gamma}$ la
normalisation de
${\mit\Gamma}$. Le morphisme $g$ se factorise en $g_1 \colon D \to
\tilde{\mit\Gamma}$ et $\varpi$. Le morphisme $g_1$
est fini et plat.
Le faisceau $g_{1*}(\mathcal{O}_D)$ est donc un fibré
vectoriel sur la $k$-courbe lisse intègre $\tilde{\mit\Gamma}$.
L'indice $n_{\tilde{\mit\Gamma}}$ de $\tilde{\mit\Gamma}$ sur $k$
est divisible par $\ell$.
En appliquant le
théorème de Riemann-Roch au fibré vectoriel $f_{1*}(\mathcal{O}_D)$
sur $\tilde{\mit\Gamma}$,
on obtient ici
$\chi_D(\mathcal{O}_D)=\chi_{\tilde{\mit\Gamma}}(g_{1*}(\mathcal{O}_D)
\equiv s\cdot \chi_{\tilde{\mit\Gamma}} (\mathcal{O}_{\tilde{\mit\Gamma}})
\pmod{\ell}$,
où $s$ est le degré du morphisme $g_1$. Comme $\chi_D(\mathcal{O}_D)$
est premier à $\ell$, on conclut qu'il en est de même
de $\chi_{\tilde{\mit\Gamma}}(\mathcal{O}_{\tilde{\mit\Gamma}})$.
On a la
suite exacte courte de faisceaux cohérents sur ${\mit\Gamma}$
suivante : $0\to \mathcal{O}_{\mit\Gamma} \to
\varpi_*\mathcal{O}_{\tilde{\mit\Gamma}}
\to \mathcal{F}\to 0$ (où le faisceau
$\mathcal{F}$ est, par définition, le
conoyau ainsi indiqué). La caractéristique d'Euler-Poincaré
du faisceau gratte-ciel $\mathcal{F}$ est la dimension du $k$-espace
vectoriel
$H^0({\mit\Gamma},
\mathcal{F})$, elle est donc divisible par $\ell$. Comme $\varpi$
est fini, on a
$\chi_{\mit\Gamma}(\varpi_*\mathcal{O}_{\tilde{\mit\Gamma}})
=\chi_{\tilde{\mit\Gamma}}(\mathcal{O}_{\tilde{\mit\Gamma}})$.
On conclut alors que $\chi_{\mit\Gamma}(\mathcal{O}_{\mit\Gamma})$
est premier à $\ell$.  Sur la surface $X$, le théorème de Riemann-Roch
pour un faisceau inversible $\mathscr{L}$, ou pour sa classe dans
$\Pic X$, s'écrit
$\chi_X(\mathscr{L})-\chi_X(\mathcal{O}_X)=
(\mathscr{L}\cdot(\mathscr{L}-\canK_X))/2$.
De la suite exacte de faisceaux cohérents $0 \to
\mathcal{O}_X(-{\mit\Gamma})
\to \mathcal{O}_X \to \mathcal{O}_{\mit\Gamma} \to 0$ sur $X$
on déduit
\[
\chi_{\mit\Gamma}(\mathcal{O}_{\mit\Gamma})=
\chi_X(\mathcal{O}_{\mit\Gamma})
= \chi_X(\mathcal{O}_X)-\chi_X(\mathcal{O}_X{{(-\mit\Gamma)}})=-
({\mit\Gamma}\cdot({\mit\Gamma}+\canK_X))/2
\]
ce qui établit le théorème.
\end{proof}

On dit qu'une surface de del Pezzo $X$ est
déployée sur le corps $L \subseteq
\bar k$ si $X(L) \neq \varnothing$
et si l'inclusion $\Pic(X_L) \hookrightarrow \Pic(\bar X)$
est une égalité.

\begin{thm}
Soient $k$ un corps de caractéristique $0$ et $X$ une
surface de del Pezzo de degré $2$ ou $4$.
Supposons que $2$ divise l'indice $n_X$ de $X$ sur $k$, que le groupe
de Picard
$\Pic X$ de $X$ est réduit aux multiples de la classe $\canK_X$ du fibré
canonique, et que la surface
$X$ est déployée par une extension galoisienne $L/k$
de degré une puissance de $2$.  Il existe alors un corps
$F$ contenant $k$, de dimension cohomologique $1$, tel que
$n_{X_F}=2$.
\label{DPKollar}
\end{thm}

\begin{proof}

Pour $K$  une extension de $k$, on note $\mathbi{P}(K)$
la conjonction des deux propriétés suivantes :
\begin{itemize}
\item[(i)] $2$ divise $n_{X_K}$,
\item[(ii)] $K\otimes_k L$ est un corps.
\end{itemize}

Soit $G$ le groupe de Galois de $K$ sur $k$.
On a les inclusions de groupes abéliens libres de rang un
$\mathbb{Z}\cdot \canK_X \hookrightarrow \Pic(X) \hookrightarrow
\Pic X_K^G \hookrightarrow
(\Pic\bar X)^{\Gal(\bar k/k)}$, et la classe de $\canK_X$
n'est pas divisible dans $\Pic\bar X$,
ces inclusions sont donc des égalités.
Sous l'hypothèse (ii), on a un isomorphisme
de $G$-modules $\Pic X_K \cong \Pic X_{K\otimes_k L}$.
On en déduit $\mathbb{Z} \cdot \canK_X= (\Pic X_K)^G = (\Pic
X_{K\otimes_k L})^G$, donc $\mathbb{Z}\cdot \canK_X = \Pic(X_L)$.

On veut maintenant appliquer le théorème \ref{InductionTheorem},
avec $p=2$.   La
condition (0) est assurée dans notre choix initial.

La propriété (i) est trivialement préservée par extension
algébrique de degré impair. La propriété (ii)
l'est aussi, de façon évidente :
ceci garantit la condition (1) du
théorème \ref{InductionTheorem}.

Supposons que $X_K$ satisfasse (i) et (ii).
Soit $C$ une
$K$-conique lisse sans point
rationnel. La caractéristique d'Euler-Poincaré
de $C$ est $1$. Soit $D$ une $K$-courbe lisse intègre,
équipée 
d'un $K$-morphisme dominant $D \to C$
de degré impair. Supposons 
$X(K(D)) \neq \varnothing$.
D'après le théorème \ref{RiemannRoch} appliqué à $\ell=2$,
l'image de la courbe $D$ dans $X_K$
est une courbe ${\mit\Gamma} \subset X_K$ telle que
$({\mit\Gamma}\cdot({\mit\Gamma}+\canK_X))/2$ est impair.
Sous l'hypothèse $\mathbb{Z} \cdot \canK_X= \Pic X_K$,
il existe un entier $n \in \mathbb{Z}$ tel que ${\mit\Gamma}=n\canK_X$,
et alors $({\mit\Gamma}\cdot({\mit\Gamma}+\canK_X))/2=  \frac{n(n+1)}{2}
(\canK_X\cdot \canK_X)$ est pair puisque
$(\canK_X\cdot \canK_X)$ est soit $4$ soit $2$.
Cette contradiction montre qu'il n'existe pas d'extension $K(D)/K(C)$ de
degré impair avec $X(K(D)) \neq \varnothing$, en d'autres termes
    $2$ divise $n_{X_{K(C)}}$.
La propriété (i) est donc préservée, en présence de la
propriété (ii), par passage au corps des fonctions d'une conique,
et la propriété (ii)
l'est aussi de façon évidente : ceci garantit la condition (2).

    Enfin, la condition (3) est claire.   Le
théorème \ref{InductionTheorem} permet alors de conclure.
\end{proof}

\textbf{Remarques :}
\par\nobreak

(1) L'existence de surfaces de del Pezzo de degré $2$ sans point
rationnel sur un corps de dimension cohomologique $1$ résulte
déjà du théorème \ref{PrincipalDP4} : on peut en effet partir
d'une surface
de del Pezzo de degré $4$ avec cette propriété
et éclater un point fermé de degré $2$ non situé
sur les $16$ droites.

\medskip
(2) Au paragraphe précédent, on a vu des exemples de
surfaces de del Pezzo de degré $4$ sur $\mathbb{Q}$ satisfaisant
les hypothèses du théorème \ref{DPKollar}.
La condition que $2$ divise $n_X$
est trivialement satisfaite pour la surface $X/\mathbb{Q}$ de
la proposition \ref{DP4Starting}, puisque cette surface
n'a pas de point sur le corps des réels.
Pour donner un exemple de surface de del Pezzo
de degré $2$, on partira d'une équation
$z^2+ax^4+by^4+cz^4=0$, avec $a, b, c \in {\mathbb{Q}}$
tous positifs et suffisamment généraux.
On renvoie au récent article \cite{KreschTschinkel}
pour le calcul de l'action du groupe de Galois
sur le groupe de Picard d'une surface de del Pezzo
de degré $2$.

\medskip

(3) Comme le remarque J. Kollár, les raisonnements
faits aux théorèmes \ref{DP4Iterating} et \ref{RiemannRoch}
ont d'autres applications. En voici une.
Soit $X$ une surface projective et lisse de degré $4$ dans
${\mathbb{P}}^3_{\mathbb{R}}$. Si $X$ est
« suffisamment générale »
et si $X$ ne possède pas de point réel, alors il n'y a aucune
courbe de genre zéro (géométriquement intègre)
tracée sur $X$.
En d'autres termes, il n'existe pas
de ${\mathbb{R}}$-morphisme d'une ${\mathbb{R}}$-conique lisse
$C$ vers $X$. De fait, si $X$ est assez générale,
un théorème de Max Noether assure que l'on a
$\Pic X= {\mathbb{Z}}\cdot H$, où $H$ est la classe d'une section
hyperplane.
Supposons qu'il existe un morphisme $C \to X$  comme évoqué.
Soit ${\mit\Gamma}$ l'image de $C$. On a alors
${\mit\Gamma}=nH \in \Pic X$. On a $\canK_X=0 \in \Pic(X)$.
Le théorème \ref{RiemannRoch} implique alors que
$2$ ne divise pas $2n^2$, ce qui est absurde : il n'existe donc
pas de tel morphisme.
On comparera ce résultat avec le résultat de
géométrie complexe assurant que sur toute surface $K3$
sur $\mathbb{C}$
il existe au moins une courbe de genre géométrique zéro.

 \vskip 2,5cm

Jean-Louis Colliot-Thélène

Mathématiques, Bâtiment 425

Université  Paris-Sud

F-91405 Orsay FRANCE

\verb=colliot@math.u-psud.fr=
    \bigskip
 
 David Madore
 
 Mathématiques, Bâtiment 425

Université  Paris-Sud

F-91405 Orsay FRANCE
 
\verb=david.madore@ens.fr=


\begin{thebibliography}{foo}
%
\bibitem{Amer}
M. Amer, \textit{Quadratische Formen über
Funktionenkörpern}, thèse (sous la direction de A. Pfister),
non publiée (1976).
%
\bibitem{Ax}
J. Ax, « A field of cohomological dimension $1$
which is not $C_1$ », \textit{Bull. Amer. Math. Soc.}
\textbf{71} (1965) 717.
%
\bibitem{Brumer}
A. Brumer, « Remarques sur les couples de formes
quadratiques »,
\textit{C. R. Acad. Sci. Paris} Sér. A-B
\textbf{286} (1978), n°16, A679--A681.
%
\bibitem{CTThmNinety}
J.-L. Colliot-Thélène, « Hilbert's Theorem 90 for
$K_2$, with
Application to the Chow Groups of Rational Surfaces »,
\textit{Invent. math.}, \textbf{71} (1983) 1--20.
%
\bibitem{CorayCub}
D. Coray, « Algebraic points on cubic
hypersurfaces », \textit{Acta
Arithmetica}, \textbf{30} (1976), 267--296.
%
\bibitem{CorayDP4}
D. Coray, « Points algébriques sur les surfaces de Del
Pezzo »,
\textit{C. R. Acad. Sci. Paris} Sér. A-B
\textbf{284} (1977), n°24, A1531--A1534.
%
\bibitem{DeJongStarr}
A. J. de Jong \& J. Starr, « Every rationally
connected variety over
the function field of a curve has a rational point »,
à paraître dans \textit{Amer. J. Math.}
%
\bibitem{Ducros}
A. Ducros, « Dimension cohomologique et points
rationnels sur les
courbes », \textit{J. Algebra}, \textbf{203} (1998)
n°2, 349--354.
%
\bibitem{Esnault}
H. Esnault, « Varieties over a finite field with
trivial Chow group of
$0$-cycles have a rational point », \textit{Invent. math.},
\textbf{151} (2003) 1, 187--191.
%
\bibitem{GraberHarrisStarr}
T. Graber, J. Harris \& J. Starr, « Families of
rationally connected
varieties », \textit{J. Amer. Math. Soc.},
\textbf{16} (2003), 57--67.
%
\bibitem{Hartshorne}
R. Hartshorne, \textit{Algebraic Geometry}, Springer, Graduate Texts in
Mathematics 52.
%
\bibitem{Kato}
K. Kato \& T. Kuzumaki, « The Dimension of Fields and
Algebraic $K$-Theory », \textit{Journal of Number Theory},
\textbf{24} (1986) 229--244.
%
\bibitem{KollarBook}
J. Kollár, \textit{Rational Curves on Algebraic Varieties}, Springer,
Ergebnisse der Mathematik und ihrer Grenzgebiete, 3. Folge, Band 32.
%
\bibitem{KreschTschinkel}
A. Kresch \& Yu. Tschinkel, « On the arithmetic
of del Pezzo surfaces of degree $2$ », prépublication
(2003).
%
\bibitem{Manin}
Yu. I. Manin, \textit{Cubic Forms: Algebra, Geometry,
Arithmetic}, North-Holland (1974, second enlarged edition 1986).
%
\bibitem{Merk-uinvar}
A. S. Merkurjev, « Simple algebras and quadratic
forms » (en russe) \textit{Izv. Akad. Nauk SSSR
Ser. Matem.},\textbf{55} (1991) ; trad. anglaise \textit{Math. USSR
Izvestiya} \textbf{38} (1992) 215--221.
%
\bibitem{MerkurjevLens}
A. S. Merkurjev, « Rost's Degree
Formula », disponible à l'adresse
\url{http://www.math.ucla.edu/~merkurev/papers/lens.dvi}
%
\bibitem{MerkurjevSteenrod}
A. S. Merkurjev, « Steenrod Operations and Degree
Formulas »,
\url{http://www.math.ucla.edu/~merkurev/papers/new.dvi}
%
\bibitem{RostDegreeOther}
A. S. Merkurjev, « Degree
Formula » (annoté par M. Rost), disponible à l'adresse
\url{http://www.math.ohio-state.edu/~rost/data/degree.pdf}
%
\bibitem{MerkurjevSuslin}
A. S. Merkurjev \& A. A. Suslin, « $K$-cohomologie des
variétés
de Severi-Brauer et homomorphisme de norme
résiduelle » (en russe),
\textit{Izv. Akad. Nauk SSSR Ser. Matem.}, \textbf{46} (1982) ;
trad. anglaise \textit{Math. USSR Izvestiya} \textbf{21} (1983)
307--340.
%
\bibitem{RostDegree}
M. Rost, « Notes on the Degree
Formula », disponible à l'adresse
\url{http://www.math.ohio-state.edu/~rost/data/bd.pdf}
%
\bibitem{SerreGACC}
J-P. Serre, \textit{Groupes algébriques et corps de classes},
Publications mathématiques de l'Université de Nancago (VII),
Actualités scientifiques et industrielles \textbf{1264}, Hermann, Paris
(1959).
%
\bibitem{SerreGC}
J-P. Serre, \textit{Cohomologie galoisienne}, Cinquième édition,
révisée
et complétée, Springer Lecture Notes in Mathematics {\bf 5} (1994)
%
\end{thebibliography}
\end{document}